\UseRawInputEncoding
\documentclass{amsart}
\usepackage{amsmath}
\usepackage{amssymb}
\usepackage{amsthm}
\usepackage{graphicx}
\usepackage{hyperref}
\usepackage[a4paper,margin=1in,footskip=0.25in]{geometry}
\newtheorem{theorem}{Theorem}

\newtheorem{lemma}[theorem]{Lemma}
\newtheorem{proposition}[theorem]{Proposition}
\newtheorem{construction}[theorem]{Construction}

\newtheorem{problem}[theorem]{Problem}
\newtheorem{remark}[theorem]{Remark}

\graphicspath{{./figures/}}
\def\d{\delta}
\def\eps{\varepsilon}

\begin{document}
	
	\title{On subgraphs of tripartite graphs}

\thanks{
This paper is a result of the Undergraduate Research Initiations in Mathematics, Mathematics Education and Statistics
(RIMMES) that the first author participated in under the guidance of the second author.
The second author is partially supported by a Simons Collaboration Grant 710094.}

\author{Abhijeet Bhalkikar}
\address{Department of Mathematics and Statistics, Georgia State University, Atlanta, GA 30303}
\email{abhalkikar1@student.gsu.edu}

\author{Yi Zhao}
\address
{Department of Mathematics and Statistics, Georgia State University, Atlanta, GA 30303}
\email{yzhao6@gsu.edu}


	\maketitle
	
	\begin{abstract}
	In 1975 Bollob\'{a}s, Erd\H{o}s, and Szemer\'{e}di \cite{1} investigated a tripartite generalization of the Zarankiewicz problem: what minimum degree forces a
    tripartite graph with $n$ vertices in each part to contain an octahedral graph $K_3(2)$?
	They proved that $n+2^{-1/2}n^{3/4}$ suffices and suggested it could be weakened to $n+cn^{1/2}$ for some constant $c>0$. In this note we show that their method only gives $n+ (1+o(1)) n^{11/12}$ and provide many constructions that show if true, $n+ c n^{1/2}$ is better possible.
	\end{abstract}
	
	\section{Introduction}
Let $K_t$ denote the complete graph on $t$ vertices.
As a foundation stone of extremal graph theory,  Tur\'{a}n's theorem in 1941 \cite{9} determines the maximum number of edges in graphs of a given order not containing $K_{t}$ as a subgraph (the $t=3$ case was proven by Mantel in 1907 \cite{7}). 
In 1975 Bollob\'{a}s, Erd\H{o}s, and Szemer\'{e}di \cite{1} investigated the following Tur\'an-type problem for multipartite graphs.
\begin{problem}
\label{pro:BES}
Given integers $n$ and $3\le t\le r$, what is the largest minimum degree $\delta(G)$ among all $r$-partite graphs $G$ with parts of size $n$ and which do not contain a copy of $K_{t}$?
\end{problem}
The $r=t$ case of Problem~\ref{pro:BES} had been a central topic in Combinatorics until it was finally settled by Haxell and Szab\'{o} \cite{MR2195582}, and Szab\'{o} and Tardos \cite{MR2246152}. Recently Lo, Treglown, and Zhao \cite{LTZ} solved many $r> t$ cases of the problem, including when $r\equiv -1 \pmod{t-1}$ and $r=\Omega(t^2)$.

For simplicity, let $G_r(n)$ denote an (arbitrary) $r$-partite graph with parts of size $n$.
Let $K_r(s)$ denote the complete $r$-partite graph with parts of size $s$. In particular, $K_3(2)$ is known as the \emph{Octahedral graph}. 
In the same paper Bollob\'{a}s, Erd\H{o}s, and Szemer\'{e}di \cite{1} also asked the following question.
\begin{problem}
\label{prob2}
Given a tripartite graph $G=G_3(n)$, what $\delta(G)$ guarantees a copy of $K_3(2)$?
\end{problem}
Problem~\ref{prob2} is a natural generalization of the well-known Zarankiewicz problem \cite{10}, whose symmetric version asks for the largest number of edges in a bipartite graph $G_2(n)$ that contains no $K_2(s)$ as a subgraph (in other words, \emph{$K_2(s)$-free}).

In \cite[Corollary 2.7]{1} the authors stated that $\d(G)\ge n + 2^{-1/2} n^{3/4}$ guarantees a copy of $K_3(2)$. This follows from \cite[Theorem 2.6]{1}, which handles the general case of $K_3(s)$ for arbitrary $s$. 
Unfortunately, there is a miscalculation in the proof of \cite[Theorem 2.6]{1} and thus 
the bound $\d(G)\ge n + 2^{-1/2} n^{3/4}$ is unverified. We follow the approach of
\cite[Theorem 2.6]{1} and obtain the following result.

%	Bollob\'{a}s, Erd\H{o}s, and Szemer\'{e}di \cite{1} proved that if $\delta(G_{3}(n)) \geq n +2^{\frac{1}{3} - \frac{1}{s^{2}}}n^{1-\frac{1}{s^{2}}}$, then $G_{3}(n)$ contains a $K_{3}(s)$. In particular, this implies that if $\delta(G_{3}(n)) \geq n+2^{-1/2}n^{3/4}$, then $G_{3}(n)$ contains a $K_{3}(2)$.

	\begin{theorem}\label{thm:ub}
		Given an integer $s\ge 2$ and $\eps > 0$, let $n$ be sufficiently large. If $G=G_3(n)$ satisfies 
		$\delta(G) \geq n+ (1+ \eps)(s-1)^{{1}/(3s^{2})}n^{1-{1}/(3s^{2})}$, then $G$ contains a copy of $K_{3}(s)$.
	\end{theorem}
	
	In particular, Theorem~\ref{thm:ub} implies that every $G= G_3(n)$ with $\delta(G)\ge n + (1+o(1)) n^{11/12}$ contains a copy of $K_{3}(2)$. Using a result of Erd\H{o}s on hypergraphs \cite{3}, we give a different proof of Theorem~\ref{thm:ub} under a slightly stronger condition $\delta(G) \geq n+ (3n)^{1- 1/(3s^2)}$. Thus $c\, n^{11/12}$ is a natural additive term for Problem~\ref{prob2} under typical approaches for extremal problems. 
	
	\medskip
	On the other hand, the authors of \cite{1} conjectured that $\delta(G)\ge n + c n^{1/2}$ suffices for Problem~\ref{prob2}. Although not explained in \cite{1}, they probably thought of Construction~\ref{cons1}, a natural construction based on the one for the Zarankiewicz problem. We indeed find many non-isomorphic constructions, Construction~\ref{cons2}, with the same minimum degree. 
	
	\begin{proposition}\label{thm:lb}
		For any $n=q^2+q+1$ where $q$ is a prime power, there are many tripartite graphs $G = G_{3}(n)$ such that $\delta(G)\ge n + n^{1/2}$ and $G$ contains no $K_{3}(2)$.
	\end{proposition}

Theorem~\ref{thm:ub} and Proposition~\ref{thm:lb} together show that the answer for Problem~\ref{prob2} lies between $n+n^{1/2}$ and $n+n^{11/12}$. The truth may be closer to the lower bound. If this is the case, 
%However, if $n+n^{1/2}$ is the correct answer, 
then verifying it may be hard given the presence of many non-isomorphic constructions.

%The best upper and lower bounds for the more general problem on $K_3(s)$ are further apart.
We know less about the minimum degree of $G_3(n)$ that forces a copy of $K_3(s)$.
Theorem~\ref{thm:ub} shows that $\d(G_3(n))\ge n + c n^{1 - 1/(3s^2)}$ suffices. 
As shown in Remark~\ref{rem}, if there is a $K_2(s)$-free bipartite graph $B= G_2(n)$ with $\d(B)=\Omega(n^{1 - 1/s})$, then our constructions for Proposition~\ref{thm:lb}
provide a tripartite $K_3(s)$-free graph $G=G_3(n)$ with $\d(G)= n + \Omega(n^{1- 1/s})$.

% Tur\'{a}n's theorem states that the complete $t$-partite graph on $n$ vertices with parts of size $\lceil \frac{n}{t} \rceil$ or $\lfloor \frac{n}{t} \rfloor$, has the most edges among all $K_{t+1}$-free graphs on $n$ vertices. 
	
%	A well known open problem in extremal graph theory is the Zarankiewicz problem \cite{10}. Let $K_{s,t}$ denote the complete bipartite graph whose one part has $s$ vertices and the other part has $t$ vertices. The Zarankiewicz problem asks for  $z(n, t)$, the largest number of edges in a bipartite graph with $n$ vertices in each part and which contains no $K_{s,t}$ as a subgraph. 

%	To obtain the lower bound for $z(n, 2)$, the bipartite graph with $n$ vertices on each side and no $K_{2,2}$ subgraph may be obtained as the Levi graph \cite{6}, which is a bipartite representation of a projective geometry configuration.
	
%	In this paper we investigate a tripartite Tur\'an problem, which extends the Zarankiewicz problem from bipartite graphs to tripartite graphs.
	
%	$G_{3}(n)$ is a tripartite graph with $3$ parts and with $n$ vertices in each part. $K_{3}(2)$ is the complete tripartite graph with $2$ vertices in each part. $\delta(G_{3}(n))$ is the minimal degree of the tripartite graph $G_{3}(n)$. Given positive integers $s\le n$, we define $f(n, s)$ to be the smallest integer $m$ such that every tripartite graph $G=G_3(n)$ with $\delta(G)> m$ contains a copy of $K_3(s)$.

	\section{Proof of Theorem~\ref{thm:ub}}
	
	In order to prove Theorem~\ref{thm:ub}, we need the following results from \cite{1}.
	\begin{lemma} \cite[Theorem 2.3]{1} \label{thm:Thm_2_3}
		Suppose every vertex of $G = G_{3}(n)$ has degree at least $n+t$ for some integer $t \leq n$. Then there are at least $t^{3}$ triangles in $G$.
	\end{lemma}
	\begin{lemma} \cite[Lemma 2.4]{1}
		\label{lemma:Lemma_2_4}
		Let $X = \{1,\ldots,N\}$ and $Y=\{1,\ldots,p\}$. Suppose 
		$A_1, \dots, A_p$ are subsets of $X$ such that
		$\sum_{i=1}^{p}|A_{i}| \geq pwN$ and $(1-\alpha)wp \geq q$, where $0<\alpha<1$ and $N$, $p$ and $q$ are natural numbers. Then there are $q$ subsets $A_{i_{1}},\ldots, A_{i_{q}}$ such that $|\bigcap\limits^q_{j=1}A_{i_{j}}| \geq N(\alpha w)^{q}$.
	\end{lemma}
	
	Let $z(n, s)$ denote the largest number of edges in a bipartite $K_2(s)$-free graph with $n$ vertices in each part.
	%and which contains no $K_2(s)$ as a subgraph. 
	K\H{o}v\'{a}ri, S\'{o}s, and Tur\'{a}n \cite{5} gave the following upper bound for $z(n, s)$.\footnote{In \cite{1} the authors instead used the Tur\'an number ex$(2n, K_2(s))$, which gives a slightly worse constant here.}
	\begin{lemma}
		\label{lem:KST}
		$z(n, s)\le (s-1)^{1/s}(n - s+1) n^{1 - 1/s} + (s-1)n$.
	\end{lemma}

	\medskip
	
	We are ready to prove Theorem~\ref{thm:ub}.
	
	\begin{proof}[Proof of Theorem~\ref{thm:ub}]
	Let $G$ be a tripartite graph with three parts
	$V_{1}, V_{2}, V_{3}$ of size $n$ each.
	Suppose $\delta(G)\ge n+ t$, where $t= (1+\eps)(s-1)^{\frac{1}{3s^{2}}}n^{1-\frac{1}{3s^{2}}} > n^{1 - \frac1{3s^2}}$. By Lemma~\ref{thm:Thm_2_3}, $G$ contains at least $t^{3}$ triangles.
	%Without loss of generality, assume that $\eps < \eps_0$ such that $1+\eps \ge (1+\eps)^{2/3}$. 

	We apply Lemma~\ref{lemma:Lemma_2_4} in the  following setting.
	Let $Y = V_{1} = \{1,\ldots,n\}$ and $X= V_2 \times V_3$ be the set of $n^{2}$ pairs $(x,y)$, $x \in V_{2}, y \in V_{3}$. For $1\le i\le n$, let $A_{i}$ be the set of pairs $(x,y) \in X$ for which $\{ i,x,y\}$ spans a triangle of $G$.
    Then  $\sum_{i=1}^n |A_i|$ is the number of triangles in $G$ so $\sum_{i=1}^n |A_i|\ge t^3$.
    Let $N=n^{2}$, $p=n$, $q=s$, $w = t^{3}/n^{3}$, and $\alpha = {1}/(1+ \eps)$. The assumptions of Lemma~\ref{lemma:Lemma_2_4} hold 
    because $pwN = t^{3}$ and 
%  $t> n^{1 - \frac1{3s^2}}$ implies that 
    \[
    (1-\alpha) w p  = \frac{\eps}{1+\eps} \left(\frac{t}{n}\right)^3 n > \frac{\eps}{1+\eps} \, n^{-1/{s^2}}\, n > s
    \]
as $n$ is sufficiently large.
   %    implies $ (1-\alpha) w p  = \frac{\eps}{1+\eps} \frac{t^3}{n^3} n  \gg s$ from \eqref{eq:t}.
By Lemma~\ref{lemma:Lemma_2_4}, there are $i_1, \dots, i_s\in V_1$ such that 
\begin{align*}
    \left|\bigcap^s_{j=1} A_{i_{j}}\right| &\ge N(\alpha w)^{q}
= n^{2}\left(\frac{t^{3}}{(1+\eps). n^{3}}\right)^{s}
%> (1+\eps)^s (s-1)^{{1}/{s}} n^{2-{1}/{s}} \quad \text{by \eqref{eq:t}}\\
%&\ge (s-1)^{{1}/{s}}n^{2-{1}/{s}} + (s-1)n 
%\quad \text{because $n$ is sufficiently large} \\
%&\ge z(n, s) \quad \text{by Lemma~\ref{lem:KST}}.
\end{align*}
Since 	
\[			t > (1+\eps)^{\frac{2}{3}}(s-1)^{\frac{1}{3s^{2}}}n^{1-\frac{1}{3s^{2}}} \quad 
			\text{and} \quad 
			\dfrac{t^{3}}{(1+\eps) n^3} > (1+\eps) (s-1)^{{1}/{s^{2}}}n^{-{1}/{s^{2}}}, 
\]
we have 
\begin{align}
	 \label{eq:B}
\left|\bigcap^s_{j=1} A_{i_{j}}\right|
> (1+\eps)^s (s-1)^{{1}/{s}} n^{2-{1}/{s}} \ge (s-1)^{{1}/{s}}n^{2-{1}/{s}} + (s-1)n. 
\end{align}

Let $B$ denote the bipartite graph between $V_2$ and $V_3$ with $E(B)= \bigcap^s_{j=1} A_{i_{j}}$. 
By \eqref{eq:B} and Lemma~\ref{lem:KST}, $B$ contains a copy of $K_2(s)$.
Since every edge of $B$ forms a triangle with each of $i_1, \dots, i_s\in V_1$, this copy of $K_2(s)$  
together with $i_1, \dots, i_s$ span a desired copy of $K_3(s)$ in $G$.
\end{proof}

	\bigskip

	We now give another proof of Theorem~\ref{thm:ub} with slightly larger $\delta(G)$ 
	by a classical result of Erd\H{o}s on hypergraphs \cite{3}. An $r$-uniform hypergraph or $r$-graph is a hypergraph such that all its edges contain exactly $r$ vertices. Let $K^r_r(s)$ denote the 
	complete $r$-partite $r$-graph with $s$ vertices in each part, namely, its vertex set consists of disjoint parts $V_1, \dots, V_r$ of size $s$, and edges set consists of all $r$-sets $\{v_1, \dots, v_r\}$ with $v_i\in V_i$ for all $i$.
	
	\begin{lemma}\cite[Theorem 1]{3} 
		\label{Erdos2}
%		Let $ex(n,F)$ denote the Turan number of the $r$-graph $F$, namely, the maximum size of an $n$-vertex $r$-graph not containing $F$ as a subgraph. Then we have that $ex(n,K^{r}(l, \ldots, l)) \leq n^{r -l^{1-r}}$.
%		Let $f(n; K^{(r)}(l_{1},\ldots,l_{r}))$ be the smallest integer such that every $G^{(r)}(n; f(n; K^{(r)}(l_{1},\ldots,l_{r})))$ contains a $K^{(r)}(l_{1},\ldots,l_{r})$. Then we have that $f(n; K^{(r)}(l_{1},\ldots,l_{r})) \leq n^{r -l^{1-r}}$.
	Given integers $r, s\ge 2$, let $n$ be sufficiently large. Then every $r$-graph on $n$ vertices with at least $n^{r -s^{1-r}}$ edges contains a copy of $K^r_r(s)$.
	\end{lemma}

	\begin{proposition}\label{prop2}
		Let $s\ge 2$ and $n$ be sufficiently large. 
		Every tripartite graph $G = G_{3}(n)$ with $\delta(G) \geq n+ (3n)^{1- 1/(3s^2)}$ contains a copy of $ K_{3}(s)$.
	\end{proposition}	
	\begin{proof}
		Suppose $G = G_{3}(n)$ satisfies $\delta(G) \geq n+ (3n)^{1- 1/(3s^2)}$. 
		By Lemma~\ref{thm:Thm_2_3}, $G$ contains at least $(3n)^{3- 1/s^2}$ triangles.
		Let $H$ be the $3$-graph on $V(G)$, whose edges are triangles of $G$. Then $H$ has $3n$ vertices and at least $(3n)^{3- s^{-2}}$ edges. By Lemma~\ref{Erdos2} with $r=3$ and $s=2$, $H$ contains a copy of $K^3_3(s)$, which gives a copy of $K_3(s)$ in $G$.
	\end{proof}
		
%	\bigskip
	
	\section{Proof of Proposition~\ref{thm:lb}}
	%Constructions for the lower bound}
	In this section we prove Proposition~\ref{thm:lb} by constructing many tripartite $K_{3}(2)$-free graphs $G_{3}(n)$ with $\delta(G_3(n))\ge n + n^{1/2}$. 
	
	One main building block is a bipartite $K_2(2)$-free graph $G_0=G_2(n)$ with  $\delta(G_0)\ge \sqrt{n}$. First shown in \cite{8}, such a graph exists when $n=q^2 + q +1$ and a projective plane of order $q$ exists. Indeed, two parts of $V(G)$ correspond to the points and lines of the projective plane and a point is adjacent to a line if and only if the point lies on the line. It is easy to see that such graph contains no $K_2(2)$ and is regular with degree $q+1> \sqrt{n}$.
	
	\begin{construction}
		\label{cons1}
		Suppose $G= G_{3}(n)$ has parts $V_{1}, V_{2}$ and $V_{3}$ each of size $n$.
		Let the bipartite graphs between $V_{1}$ and $V_{2}$ and between $V_{1}$ and $V_{3}$ be complete, while the bipartite graph between $V_{2}$ and $V_{3}$ is $G_0$ defined above.
		
		Since $\deg_{G_0}(v) \ge \sqrt{n}$ for $v \in V_{2} \cup V_{3}$,
		we have $\delta(G) \ge n + \sqrt{n}$. Furthermore, 
		$G$ contains no $K_{3}(2)$ because by the definition of $G_0$, there is no $K_{2}(2)$ between $V_{2}$ and $V_{3}$.
	\end{construction}
	
	\bigskip
	
	We now provide a family of constructions with the same properties.
	\begin{construction}\label{cons2}
		\begin{figure}[hbt!]
			\centering
			\includegraphics[scale=0.25]{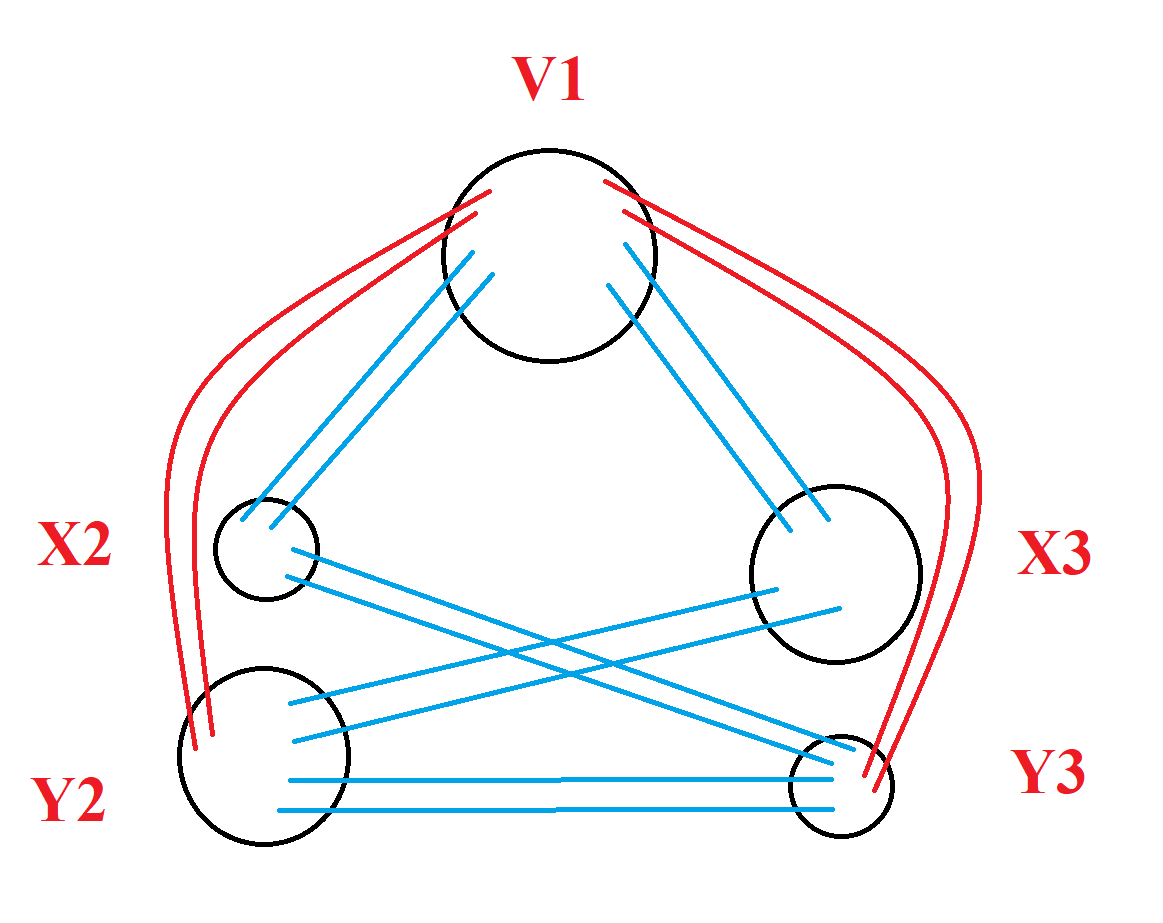}
			\caption{Graph from Construction~\ref{cons2}}
			\label{fig:construction2}
		\end{figure}
		Let $G=G_3(n)$ be a tripartite graph with parts $V_1$, $V_2$, and $V_3$ of size $n$ each. Partition $V_2= X_2\cup Y_2$ arbitrarily such that $\sqrt{n}\le |X_2|\le |Y_2|$.  Partition $V_3= X_3\cup Y_3$ arbitrarily such that $|X_3|= |Y_2|$ and $|Y_3|= |X_2|$. 
		
		The bipartite graphs $(V_1, X_2)$, $(X_2, Y_3)$, $(Y_3, Y_2)$, $(Y_2, X_3)$, and $(X_3, V_1)$ are complete, in other words, $V_1, X_2, Y_3, Y_2, X_3$ form a blowup of $C_5$. Let the bipartite graph between $V_1$ and $Y_2\cup Y_3$ be isomorphic to $G_0$ (note that $|X_2|+|Y_2|=|X_3|+ |Y_3|= n$).
		
		For any vertex $v\in X_2$, $\deg(v) = |V_1|+|Y_3|\ge n+ \sqrt{n}$. The vertices $v\in X_3$ satisfy $\deg(v) = |V_1|+|Y_2|\ge n+ n/2$. For any $v\in Y_2$, $\deg(v)\ge |V_3|+ \delta(G_0) \ge n + \sqrt{n}$. The same holds for the vertices of $Y_3$. At last, every vertex $v\in V_1$ satisfies $\deg(v)\ge |X_2|+|X_3|+ \delta(G_0)\ge n+ \sqrt{n}$. These together show that $\delta(G)\ge n+ \sqrt{n}$.
		
		Suppose $G$ contains a copy of $K_3(2)$ with vertex set $S$.
		Then $|S\cap V_i|= 2$ for $i=1, 2, 3$. Since there is no edge between $X_2$ and $X_3$, either $S\cap X_2= \emptyset$ or $S\cap X_3 = \emptyset$. Suppose, say, $S\cap X_2= \emptyset$, which forces $|S\cap Y_2|= 2$. Hence $S\cap Y_2$ and $S\cap V_1$ span a copy of $K_2(2)$, contradicting the definition of $G_0$. 
		
	\end{construction}
	
	If letting $X_2=\emptyset = Y_3$ in Construction~\ref{cons2}, then we obtain Construction~\ref{cons1}. Nevertheless, we prefer viewing 
	Constructions~\ref{cons1} and \ref{cons2} as different constructions because after removing $o(n^2)$ edges, Construction~\ref{cons2} contains many 5-cycles while Construction~\ref{cons1} does not.
	
	\begin{remark}\label{rem}
	    If we replace $G_0$ by a $K_2(s)$-free bipartite graph with $n$ vertices in each part in Constructions~\ref{cons1} and \ref{cons2}, then we obtain a $K_3(s)$-free tripartite graph $G_3(n)$. It has been conjectured that
	    there exist a $K_2(s)$-free bipartite graph with $n$ vertices in each part and $\Omega( n^{2 - 1/s})$ edges (this is known for $s=2, 3$ \cite{2, 8}). If such bipartite graph exists and is regular, then (revised) Constructions~\ref{cons1} and \ref{cons2} provide a $K_3(s)$-free tripartite graph $G= G_3(n)$ with $\d(G)= n + \Omega( n^{1- 1/s} )$.
	\end{remark}

	%\newpage

\end{document}